\documentclass[a4paper,12pt]{amsart}

\usepackage{amscd}
\usepackage{amsmath}

\usepackage{amsthm}

\usepackage{graphicx}

\usepackage{txfonts}
\usepackage{dsfont}
\usepackage{mathabx}

\usepackage{tikz-cd}
\usepackage[titletoc]{appendix}

\usepackage{geometry}
\usepackage{hyperref}
\usepackage{tikz}

\theoremstyle{definition}
\newtheorem{Defn}{Definition}[section]

\theoremstyle{plain}
\newtheorem{Lemma}[Defn]{Lemma}

\newtheorem{theorem}[Defn]{Theorem}
\newtheorem{Corollary}[Defn]{Corollary}
\newtheorem*{main}{Main Theorem}

\theoremstyle{remark}
\newtheorem{remark}[Defn]{Remark}
\newtheorem{Example}[Defn]{Example}

\title{On a pair of extensions of Mayer-Vietoris functors}
\author{Minkyu Kim}
\date{}

\begin{document}

\maketitle
\begin{abstract}
The purpose of this paper is to construct {\it a pair} of some canonical extensions of Mayer-Vietoris functors (e.g. homology or cohomology theories of spaces).
We derive such extensions by revealing a partial universality of a cospan category of CW-spaces which is a homotopy-theoretic analogue of cobordism categories.
Those extensions are closely related with the classical theory part of abelian Dijkgraaf-Witten-Freed-Quinn model and bicommutative Turaev-Viro-Barrett-Westbury model.
\end{abstract}
\tableofcontents

%%%%%%%%%%%%%%%%%%%%%%%%%%%%%%%%

\section{Introduction}

The purpose of this paper is to construct {\it a pair} of some canonical extensions of Mayer-Vietoris functors (e.g. homology or cohomology theories of spaces).
We derive such extensions by revealing a partial universality of a cospan category of CW-spaces which is a homotopy-theoretic analogue of cobordism categories.
Those extensions are closely related with abelian Dijkgraaf-Witten-Freed-Quinn model \cite{DW}\cite{FQ}\cite{freed2009topological}\cite{monnier2015higher} and bicommutative Turaev-Viro-Barrett-Westbury model \cite{TV}\cite{BW}\cite{barrett1995equality}.
In particular, our main theorem completes some missing proofs in subsection 4.2 of our previous paper\footnote{There is a remark about terminologies. See Remark \ref{202205261256}} \cite{kim2021pair}  where we refer a normalization of DWFQ and TVBW invariants to an obstruction problem.

We explain some notations for our main theorem.
Let $\mathcal{A}$ be an abelian category and denote by $\mathsf{Rel} ( \mathcal{A})$ the relation category associated with $\mathcal{A}$ \cite{mac1961algebra} \cite{puppe1962korrespondenzen} \cite{hilton1966correspondences}.
For an extended natural number $N = 0, 1, \cdots , \infty$, we consider a category $\mathsf{Cosp}^\simeq_{\leq N} ( \mathsf{CW}^\mathsf{fin}_\ast)$ generated by $r$-dimensional cospans of finite CW-spaces up to a homotopy for $r \leq N$ (see Definition \ref{202201301536}).

\begin{main}
Let $E$ be an $\mathcal{A}$-valued Mayer-Vietoris functor for pointed finite CW-spaces with dimension lower than or equal to $(d+1)$.
We obtain a pair of canonical functors $\hat{E}$, $\check{E}$ from $\mathsf{Cosp}^\simeq_{\leq (d+1)} ( \mathsf{CW}^\mathsf{fin}_\ast)$ to $\mathsf{Rel} ( \mathcal{A})$ which satisfy (1) and (2) respectively.
Furthermore, we have a coincidence property in (3).
\begin{enumerate}
\item
The functor $\hat{E}$ is a {\it unique} dagger-preserving functor with the following commutative diagram.
Here, $i_d^\ast E$ is the restriction of $E$ to pointed finite CW-spaces with dimension $\leq d$.
The functor $\hat{E}$ is called the cospanical extension of $i_d^\ast E$.
\begin{equation}
\notag
\begin{tikzcd}[sep=scriptsize]
\mathsf{Ho} \left( \mathsf{CW^{fin}}_{\ast, \leq d} \right) \ar[d, hookrightarrow]  \ar[r, "i_d^\ast E"] & \mathcal{A} \ar[d, hookrightarrow] \\
\mathsf{Cosp}^\simeq_{\leq (d+1)} ( \mathsf{CW}^\mathsf{fin}_\ast ) \ar[r, "\hat{E}"] & \mathsf{Rel} ( \mathcal{A})
\end{tikzcd}
\end{equation}
\item
The functor $\check{E}$ is a {\it unique} dagger-preserving functor with the following commutative diagram.
Here, $\Sigma_d^\ast E$ is the pullback of $E$ by the suspension.
The functor $\check{E}$ is called the spanical extension of $\Sigma_d^\ast E$.
\begin{equation}
\notag
\begin{tikzcd}[sep=scriptsize]
\mathsf{Ho} \left( \mathsf{CW^{fin}}_{\ast, \leq d} \right) \ar[d, hookrightarrow]  \ar[r, "\Sigma_d^\ast E"] & \mathcal{A} \ar[d, hookrightarrow] \\
\mathsf{Cosp}^\simeq_{\leq (d+1)} ( \mathsf{CW}^\mathsf{fin}_\ast ) \ar[r, "\check{E}"] & \mathsf{Rel} ( \mathcal{A})
\end{tikzcd}
\end{equation}
\item
If $d \geq 1$, then the cospanical extension $\hat{E^\prime}$ extends to the spanical extension $\check{E}$ where $E^\prime = \Sigma_d^\ast E$ :
\begin{equation}
\notag
\begin{tikzcd}[sep=scriptsize]
\mathsf{Cosp}^\simeq_{\leq d} ( \mathsf{CW}^\mathsf{fin}_\ast ) \ar[r, "\hat{E^\prime}"] \ar[d, hookrightarrow] & \mathsf{Rel} ( \mathcal{A}) \ar[d, equal] \\
\mathsf{Cosp}^\simeq_{\leq (d+1)} ( \mathsf{CW}^\mathsf{fin}_\ast ) \ar[r, "\check{E}"] & \mathsf{Rel} ( \mathcal{A})
\end{tikzcd}
\end{equation}
\end{enumerate}
\end{main}

The proof follows from Corollary \ref{202201271352}.
Our strategy is to study a characterization of cospan categories of CW-spaces.
In this paper, we give a partial characterization relatively to the homotopy category of CW-spaces (see Theorem \ref{202201172145}).

We have a symmetric monoidal version of the above results.
In fact, $E$ is pointed, i.e. $E(\ast )= 0$, if and only if the (co)spanical extensions are canonically symmetric monoidal.
See Remark \ref{202201301549}.

The main theorem provides a pair of systematic ways to obtain TQFT or its homotopy-theoretic analogue from Mayer-Vietoris functors.
For a symmetric monoidal category $\mathcal{C}$, a symmetric monoidal {\it projective} functor $F :\mathsf{Rel} ( \mathcal{A}) \to \mathcal{C}$ induces a {\it $\mathcal{C}$-valued relative TQFT} or a {\it $\mathcal{C}$-valued TQFT with nontrivial anomaly} by a composition with the (co)spanical extension.
The normalization of such TQFT's or the trivialization of anomaly derives an algebraic problem which we solve in \cite{kim2021pair} for a certain $F$ based on the coincidence property.

The cospanical extension (the projective functor $F$, resp.) implicitly and partially appears in the form of {\it classical field theory} ({\it finite path integral}, resp.) in the literature \cite{FQ}\cite{freed2009topological}\cite{grandis2007collared}\cite{morton2015cohomological}\cite{sharma2017categorification}.
On the other hand, the spanical extension reveals a new way of finite path-integrals along a {\it span}  induced by a cobordism  a cospan of spaces (see Definition \ref{201912241125}).

The cospan category of CW-spaces in this paper could be regarded as a variation of the collarable cospan category \cite{grandis2007collared}.

The relation category $\mathsf{Rel} (\mathcal{A} )$ naturally appears in homological algebra, and it was systematically studied based on general categories.
The category was introduced by relations in $\mathcal{A}$ and characterized by some axioms \cite{mac1961algebra}\cite{puppe1962korrespondenzen} ; or it could be obtained from spans in $\mathcal{A}$ \cite{hilton1966correspondences}.
Furthermore, the construction is generalized to relations in any regular category \cite{meisen1974relations}.
In this paper, we follow the framework by \cite{hilton1966correspondences} which is a technically different setting.

We present a self-contained description of the relation category $\mathsf{Rel} (\mathcal{A} )$ and the cospan category $\mathsf{Cosp}^\simeq_{\leq d} ( \mathsf{CW}^\mathsf{fin}_\ast )$ in this paper.

There is a similar result to the first part of our main theorem in section 4 of \cite{grandis2007collared} based on $\mathcal{A}= \mathsf{Ab}$ and a (generalized) homology theory $E$.
Apart from slight differences in setting, Grandis deals with double categories of cospans of collarable spaces and relations in $\mathsf{Ab}$ whereas we consider only ordinary categories.
We believe that this paper could be improved to such double categories and functors in an obvious way.

\subsection{Organization}

This paper is organized as follows.
In subsection \ref{202201301622}, we give a definition of cospan categories of (pointed) spaces.
In subsection \ref{2022013016221}, we introduce cospan categories of CW-spaces as subcategories of the above one.
In subsection \ref{202201180923}, we give some algebraic characterization of $0,1$-dimensional cospan categories.
In section \ref{202201172154}, we define Mayer-Vietoris positive-and-negative diagrams.
In subsection \ref{202201301625}, we introduce a universality of MV PN diagrams and prove that the cospan category of pointed finite CW-spaces satisfies the universality.
In subsection \ref{202201301626}, we introduce a unification property of MV PN diagrams with finite dimensions and prove that finite-dimensional cospan categories of pointed finite CW-spaces satisfy the property relative to some functors.
In subsection \ref{202201301628}, we give a definition of $\mathcal{A}$-valued Mayer-Vietoris functors with some examples.
In subsection \ref{202201301328}, we give an overview of relation category of abelian categories.
In subsection \ref{202201301629}, we give a proof of our main theorem.

\section*{Acknowledgements}
The author was supported by FMSP, a JSPS Program for Leading Graduate Schools in the University of Tokyo, and JPSJ Grant-in-Aid for Scientific Research on Innovative Areas Grant Number JP17H06461.

\section{Cospans of spaces}
\label{202201141712}

In this section, we introduce some homotopy-theoretic notions related with cospans of spaces.
Especially, they induce {\it cospan categories of spaces} which are analogues of cobordism categories.
Furthermore, we explain not only a similarity but also canonical functors between cospan categories and cobordism categories (see Example \ref{202202021509}).

A cospan in a category $\mathcal{C}$ consists of three objects $X,Y,Z$ and two morphisms $f : X \to Z$ and $g : Y \to Z$.
It is usually depicted by a diagram below :
\begin{equation}
\notag
\begin{tikzcd}
X \stackrel{f}{\to} Z \stackrel{g}{\leftarrow} Y .
\end{tikzcd}
\end{equation}
In this paper, we call $Z$ by {\it bulk} and $f,g$ by {\it boundaries} of the cospan.
The category $\mathcal{C}$ of our main interest is the category of (pointed) spaces and (pointed) continuous maps.
Denote by $\mathsf{Top}_\ast$ and $\mathsf{Top}$ the categories respectively.

\begin{remark}
We remark that a {\it category} in this paper means a locally small category.
In other words, the morphisms between given two objects form a set whereas the objects could form a proper class in a naive sense.
It is only a technical assumption for our discussion to be based on ZFC.
The readers are referred to Remark \ref{202202031402}, \ref{202202031403}.
\end{remark}

\subsection{Cospan categories of spaces}
\label{202201301622}

Let $L, L^\prime$ be pointed spaces.
They are homotopy equivalent with each other if there exists pointed maps $h : L \to L^\prime$ and $h^\prime : L^\prime \to L$ such that $h^\prime \circ h \simeq \mathrm{id}_{L}$ and $h \circ h^\prime \simeq \mathrm{id}_{L^\prime}$.
On the one hand, any pointed space $L$ gives rise to a cospan in the homotopy category whose boundaries are trivial, which is described by a diagram $( \ast \to L \leftarrow \ast )$.
Here $\ast$ is the pointed one-point space.
Then the homotopy equivalence of spaces extends to an equivalence relation on cospans of spaces :

\begin{Defn}
Let $\Lambda_0, \Lambda_1$ be cospans in $\mathsf{Top}_\ast$.
The cospans $\Lambda_0, \Lambda_1$ are {\it homotopy equivalent} if there exists a homotopy equivalence between the bulks which intertwines the boundaries.
To be precise, let $\Lambda_0 = ( K_0 \stackrel{f_0}{\to} L \stackrel{f_1}{\leftarrow} K_1 )$ and  $\Lambda_1 = ( K_0^\prime \stackrel{g_0}{\to} L^\prime \stackrel{g_1}{\leftarrow} K_1^\prime)$ be cospans in $\mathsf{Top}_\ast$.
The cospans $\Lambda_0, \Lambda_1$ are {\it homotopy equivalent} if $K_0 = K_0^\prime$, $K_1 = K_1^\prime$ and there exists a homotopy equivalence $h : L \to L^\prime$ such that $g_i \simeq h \circ f_i$ for $i = 0,1$.
\end{Defn}

There is a natural binary operation on pointed spaces called the {\it wedge sum of spaces}.
It has two generalizations to cospans of spaces.
One is the wedge sum of cospans, and the other one is the homotopy pushout composition of cospans :

\begin{Defn}
Let $\Lambda_0 = ( K_0 \stackrel{f_0}{\to} L \stackrel{f_1}{\leftarrow} K_1 )$ and $\Lambda_1 = ( K_0^\prime \stackrel{g_0}{\to} L^\prime \stackrel{g_1}{\leftarrow} K_1^\prime)$ be cospans in $\mathsf{Top}_\ast$.
\begin{itemize}
\item
We define the {\it wedge sum of cospans} by the following cospan :
\begin{align}
\notag
\Lambda_0 \vee \Lambda_1
\stackrel{\mathrm{def.}}{=} ( K_0 \vee K_0^\prime \stackrel{f_0 \vee g_0}{\to} L \vee L^\prime \stackrel{f_1 \vee g_1}{\leftarrow} K_1 \vee K_1^\prime ) .
\end{align}
\item
The cospans $\Lambda_1 , \Lambda_0$ are {\it composable} if $K_1 = K_0^\prime$.
We define the {\it homotopy pushout composition} by the following cospan : 
\begin{align}
\notag
\Lambda_1 \circ \Lambda_0 \stackrel{\mathrm{def.}}{=} 
( K_0 \stackrel{h_0 \circ f_0}{\to} \mathrm{Cyl} ( f_1 , g_0 ) \stackrel{h_1 \circ g_1}{\leftarrow} K_1^\prime ) .
\end{align}
\end{itemize}
Here, the space $\mathrm{Cyl} ( f_1 , g_0 )$ is the double mapping cylinder of $f_1$ and $g_0$, and $h_0 : L \to \mathrm{Cyl} ( f_1 , g_0 )$ and $h_1 : L^\prime \to \mathrm{Cyl} ( f_1 , g_0 )$ are the canonical inclusions.
\end{Defn}

\begin{remark}
For a pointed map $f : K \to L$ between pointed spaces, we define the {\it mapping cylinder} $\mathrm{Cyl} ( f )$ by a pointed adjunction space,
\begin{align}
\notag
\mathrm{Cyl} ( f ) \stackrel{\mathrm{def.}}{=} \mathrm{Cyl} ( K ) \vee_f L .
\end{align}
In other words, it is the quotient pointed space by identifying $[x,1] \in \mathrm{Cyl} ( K) = K \wedge [0,1]^+$ with $f(x) \in L$.
By the canonical inclusion $L \to \mathrm{Cyl} ( f )$, consider $L$ as a subspace of $\mathrm{Cyl} ( f )$.
We also identify $K$ with $K \wedge \{0 \}^+ \subset K \wedge [0,1]^+ = \mathrm{Cyl} ( K)$.
Let $f : K \to L_0$ and $g : K \to L_1$ be pointed maps.
The {\it double mapping cylinder of $f$ and $g$}, denoted by $\mathrm{Cyl}(f,g)$, is defined by a pointed adjunction space $\mathrm{Cyl} ( f ) \vee_K \mathrm{Cyl} ( g )$.
It is well-known as a concrete model of homotopy pushout.
\end{remark}

\begin{remark}
The wedge sum of cospans and the homotopy pushout composition respects the homotopy equivalence of cospans of spaces.
It follows from the homotopy invariance of the wedge sum of spaces and the double mapping cylinder respetively.
\end{remark}

\begin{remark}
These notions give a generalization of the wedge sum of spaces in the sense that $t ( L_0 \vee L_1 ) = t (L_0) \vee t (L_1)$ and $t ( L_0 \vee L_1 ) = t (L_0) \circ t (L_1)$.
Here, $t (L)$ is the homotopy equivalence class of the induced cospan $( \ast \to L \leftarrow \ast )$ for an arbitrary pointed space $L$.
\end{remark}

\begin{Defn}
\label{202201271145}
We define a category $\mathsf{Cosp}^{\simeq} ( \mathsf{Top}_\ast )$ starting from cospans in the homotopy category of pointed spaces.
The objects are pointed spaces and the morphisms are homotopy equivalence classes of cospans of pointed spaces.
The identity on an object $K$ is given by the homotopy equivalence class of $( K \stackrel{\mathrm{id}}{\to} K \stackrel{\mathrm{id}}{\leftarrow} K)$.
The composition is given by the homotopy pushout composition.
Similarly, one can define a category $\mathsf{Cosp}^{\simeq} ( \mathsf{Top} )$ starting from cospans of spaces (without basepoint).
The categories $\mathsf{Cosp}^{\simeq} ( \mathsf{Top}_\ast )$ and $\mathsf{Cosp}^{\simeq} ( \mathsf{Top} )$ are naturally endowed with a dagger symmetric monoidal category structure.
We sketch the structure here.
The wedge sum of cospans induces a monoidal structure.
The usual symmetry of the wedge sum $L \vee L^\prime \cong L^\prime \vee L$ induces a symmetric monoidal structure.
Furthermore, the assignment of a transposition $( K_1 \stackrel{f_1}{\to} L \stackrel{f_0}{\leftarrow} K_0 )$ to a cospan $( K_0 \stackrel{f_0}{\to} L \stackrel{f_1}{\leftarrow} K_1 )$ provides a dagger structure.
\end{Defn}

\begin{remark}
\label{202202031402}
We give a technical remark about the size issue in set theory.
Note that the categories $\mathsf{Cosp}^\simeq ( \mathsf{Top}_\ast )$ and $\mathsf{Cosp}^\simeq ( \mathsf{Top} )$ are not defined as an mathematical object, or they are {\it navie} categories, if one only assumes ZFC.
Those categoires are not even locally small since the morphism class between two given objects is not a set.
We deal with cospan categories of (pointed) finite CW-spaces (which appears in the next subsection) throughout this paper to keep our discussion in ZFC.
Indeed, it is essentially small, in particular locally small, since the class of finite CW-spaces is essentially a set.
On the other hand, if one assumes some proper class-set theory such as NBG or the Grothendieck universe in ZFC, then the results in this paper are appropriately generalized to cospan categories of arbitrary spaces or CW-spaces.
\end{remark}

\subsection{Cospan categories of CW-spaces}
\label{2022013016221}

We introduce cospan categories of CW-spaces forming the algebraic base of this paper.
These are subcategories of $\mathsf{Cosp}^{\simeq} ( \mathsf{Top}_\ast )$ and $\mathsf{Cosp}^{\simeq} ( \mathsf{Top} )$ and have a filtration by dimension.
We first give an overview of CW-spaces and their dimension.

\begin{Defn}
A topological space $K$ (without a basepoint) is a {\it CW-space} if there exists a CW-complex structure on $K$.
Note that we do not fix any distinguished CW-complex structure.
A CW-space $K$ is {\it finite} if it has a finite CW-complex structure.
Denote by $\mathsf{CW}$ ($\mathsf{CW^{fin}}$, resp.) the categories of (finite, resp.) CW-spaces and continuous maps.
A {\it pointed CW-space} is a CW-space with a basepoint.
Denote by $\mathsf{CW}_\ast$ ($\mathsf{CW^{fin}}_\ast$, resp.) the category of (finite, resp.) CW-spaces and pointed continuous maps.
\end{Defn}

\begin{remark}
It is well-known that a CW-complex is compact if and only if its CW-complex structure is finite.
Hence, a CW-space is finite if and only if it is compact.
\end{remark}

\begin{remark}
\label{202109300008}
A typical example related with Dijkgraaf-Witten-Freed-Quinn model is given by compact smooth manifolds.
It is technically convenient for us to deal with CW-spaces instead of CW-complexes when we deal with smooth manifolds.
Indeed, a smooth manifold has a CW-complex structure which is not canonical.
The underlying space of a (compact) smooth manifold is a (finite) CW-space.
\end{remark}

\begin{remark}
In spite of the convenience, it is essentially the same with CW-complexes in the sense of homotopy.
Let $\mathbb{CW}$ be the category of CW-complexes and cellular maps.
Then the forgetful functor between the homotopy categories $\mathsf{Ho} ( \mathbb{CW} ) \to \mathsf{Ho} ( \mathsf{CW} )$ induces an equivalence of categories (similarly for pointed cases).
It is deduced from the cellular approximation theorem.
\end{remark}

\begin{Defn}
Recall that the (supremum) dimension of a CW-complex is defined by the supremum of the dimensions of cells.
We define the {\it (supremum) dimension} of a CW-space $K$ by the dimension of $X$ for a CW-complex structure on $K$.
The dimension of an empty space is defined by $(-1)$.
It is independent of the choice of $X$ due to the invariance of domain.
For an extended natural number $d = 0, 1, \cdots , \infty$, denote by $\mathsf{CW}_{\leq d}$ ($\mathsf{CW^{fin}}_{\leq d}$, resp.) the category of (finite, resp.) CW-spaces with dimension lower than or equal to $d$.
For a pointed CW-space $K$, we define the dimension by that of the underlying CW-space if $K$ does not consist of only the basepoint.
Otherwise, the dimension is defined by $(-1)$.
Denote by $\mathsf{CW}_{\ast , \leq d}$ ($\mathsf{CW^{fin}}_{\ast, \leq d}$, resp.) the category of (finite, resp.) pointed CW-spaces with dimension lower than or equal to $d$.
We make use of the symbol $\dim K$ to denote the dimension of a (pointed) CW-space $K$.
\end{Defn}

\begin{Defn}
Let $\Lambda = ( K_0 \stackrel{f_0}{\to} L \stackrel{f_1}{\leftarrow} K_1 )$ be a cospan in the category of (pointed) CW-spaces.
We define the {\it dimension} of $\Lambda$ by the supremum of $(1 + \dim K_0 )$, $(1 + \dim K_1 )$ and $\dim L$.
We introduce the symbol $\dim \Lambda$ to denote the dimension of a cospan $\Lambda$.
Similarly, we define the {\it component-dimension $\dim_\mathrm{c} \Lambda$} of $\Lambda$ by the supremum of $\dim K_0$, $\dim K_1$ and $\dim L$.
\end{Defn}

\begin{Defn}
\label{202201301536}
Let $d$ be an extended natural number, i.e. $d = 0, 1, \cdots , \infty$.
We define a subcategory $\mathsf{Cosp}^{\simeq}_{\leq d} ( \mathsf{CW}^\mathsf{fin}_\ast )$ of $\mathsf{Cosp}^{\simeq} ( \mathsf{Top}_\ast )$.
The objects are pointed finite CW-spaces with dimension lower than or equal to $(d -1)$.
The morphisms are homotopy equivalence classes of cospans $\Lambda$ with $\dim \Lambda \leq d$.
It follows from definitions that it is closed under the composition.
In a parallel way, cospans of finite CW-spaces (without basepoint) define a subcategory $\mathsf{Cosp}^{\simeq}_{\leq d} ( \mathsf{CW}^\mathsf{fin} )$ of $\mathsf{Cosp}^{\simeq} ( \mathsf{Top} )$.
Note that these subcategories inherit a dagger symmetric monoidal category structure.
\end{Defn}

\begin{Example}
\label{202202021509}
We relate cospan categories with cobordism categories when one deals with Dijkgraaf-Witten-Freed-Quinn model and Turaev-Viro-Barrett-Westbury model.
If $n$ is a natural number such that $n \leq d$, then an $n$-dimensional smooth cobordism induces a $d$-dimensional cospan of finite CW-spaces by Remark \ref{202109300008}.
It gives a rise to a (dagger symmetric monoidal) functor $\mathsf{Cob}_n \to \mathsf{Cosp}^{\simeq}_{\leq d} ( \mathsf{CW}^\mathsf{fin} )$.
\end{Example}

\subsection{Low-dimensional cospan categories}
\label{202201180923}

For convenience, we give some algebraic characterization of low-dimensional cospan categories of CW-spaces here.
Our strategy is based on generator-relation description.

\subsubsection{0-dimension}

The assignment of $K^+ = K \amalg \ast$, i.e. the disjoint union with a formal point, to a space $K$ without a basepoint induces a functor $\mathsf{Cosp}^{\simeq}_{\leq 0} ( \mathsf{CW}^\mathsf{fin} ) \to \mathsf{Cosp}^{\simeq}_{\leq 0} ( \mathsf{CW}^\mathsf{fin}_\ast )$.
It turns out to give a category equivalence.
The category $\mathsf{Cosp}^{\simeq}_{\leq 0} ( \mathsf{CW}^\mathsf{fin} )$ is algebraically described as follows.
We define a category $\bullet // \mathbb{N}$ by the category induced by the monoid structure of natural numbers which is unpacked as follows.
It has only one object $\bullet$ and morphisms are given by natural numbers $n : \bullet \to \bullet$.
The composition is given by $n \circ m  = (n+m)$ and the identity is $0 : \bullet \to \bullet$.
Then the assignment of $| L |$ (the order of $L$) to a cospan $\left( \phi \to L \leftarrow \phi \right)$ having dimension lower than or equal to $0$ induces a category equivalence $\mathsf{Cosp}^{\simeq}_{\leq 0} ( \mathsf{CW}^\mathsf{fin} ) \to \bullet // \mathbb{N}$.

\subsubsection{1-dimension}

\begin{figure}[h]
  \includegraphics[width=6.7cm]{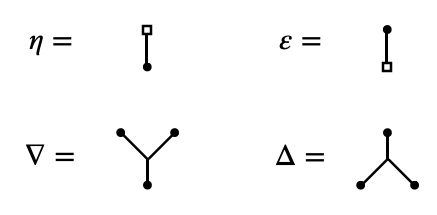}
  \caption{Graphical description of some elementary cospans with $\dim \leq 1$}
  \label{202111281133}
\end{figure}

\begin{figure}[h]
  \includegraphics[width=8cm]{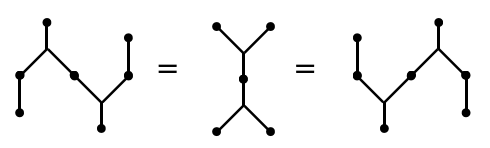}
  \caption{Frobenius relations}
  \label{2020111911432}
\end{figure}

Let $\mathrm{pt}$ be a one-point space.
In the category $\mathsf{Cosp}^{\simeq}_{\leq 1} ( \mathsf{CW}^\mathsf{fin} )$, there are some elementary morphisms $\eta, \nabla$ induced by the trivial map $\phi \to \mathrm{pt}$ and the collapsing map $\mathrm{pt} \amalg \mathrm{pt} \to \mathrm{pt}$ respectively.
Here the cospan induced by a map $f : K \to L$ is defined by $( K \stackrel{f}{\to} L \stackrel{\mathrm{id}}{\leftarrow} L)$.
Their daggers give morphisms $\varepsilon = \eta^\dagger : \mathrm{pt} \to \phi$ and $\Delta =\nabla^\dagger : \mathrm{pt} \to \mathrm{pt} \amalg \mathrm{pt}$.
We depict these morphisms as diagrams in Figure \ref{202111281133} where the black dot is the object $\mathrm{pt}$ and the blank dot is the empty space $\phi$.
By some homotopy equivalences, one could obtain a commutative Frobenius object with a unit $\eta$, a counit $\varepsilon$, a multiplication $\nabla$ and a comultiplication $\Delta$.
For example, the Frobenius relation is shown by Figure \ref{2020111911432}.
In fact, it turns out that $\mathsf{Cosp}^{\simeq}_{\leq 1} ( \mathsf{CW}^\mathsf{fin} )$ is a symmetric monoidal category freely generated by a commutative Frobenius object.
The proof is parallel with that for the 2-dimensional oriented cobordism category \cite{kock2004frobenius}.

\begin{figure}[h]
  \includegraphics[width=6cm]{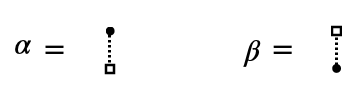}
  \caption{Graphical description of some elementary cospans with $\dim \leq 1$}
  \label{202201142148}
\end{figure}

\begin{figure}[h]
  \includegraphics[width=8cm]{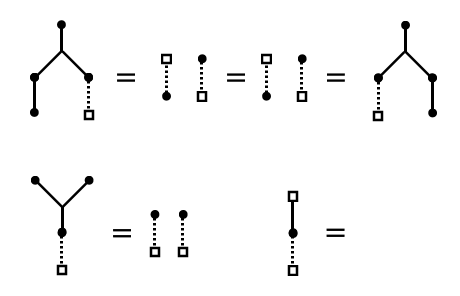}
  \caption{Relations related with $\alpha, \beta$}
  \label{202201142147}
\end{figure}

In the 1-dimensional case (and higher dimensional cases), there is a difference between pointed spaces and spaces without basepoint.
Indeed, the functor $\mathsf{Cosp}^{\simeq}_{\leq 1} ( \mathsf{CW}^\mathsf{fin} ) \to \mathsf{Cosp}^{\simeq}_{\leq 1} ( \mathsf{CW}^\mathsf{fin}_\ast )$ which adds a formal point is essentially surjective and faithful but not full.
Let $\eta : \ast \to \mathrm{pt}^+, \varepsilon : \mathrm{pt}^+ \to \ast, \nabla : \mathrm{pt}^+ \vee \mathrm{pt}^+ \to \mathrm{pt}^+$ and $\Delta : \mathrm{pt}^+ \to \mathrm{pt}^+ \vee \mathrm{pt}^+$ be the induced morphisms by the above corresponding morphisms with an abuse of notations.
Here, $\ast$ is a pointed one-point space.
We set one more morphism $\alpha : \mathrm{pt}^+ \to \ast$ to be the morphism in $\mathsf{Cosp}^{\simeq}_{\leq 1} ( \mathsf{CW}^\mathsf{fin}_\ast )$ induced by the collapsing map $\mathrm{pt}^+ \to \ast$.
It is depicted in Figure \ref{202201142148}.
It is easy to verify the relations in Figure \ref{202201142147} where the empty graph in the last relation means the identity on $\ast$.
It turns out that $\mathsf{Cosp}^{\simeq}_{\leq 1} ( \mathsf{CW}^\mathsf{fin}_\ast )$ is a symmetric monoidal category freely generated by a commutative Frobenius object equipped with $\alpha , \beta = \alpha^\dagger$ subject to the relations in Figure \ref{202201142147}.

%%%%%%%%%%%%%%%%%%%%%%
\section{Characterization problem of cospan categories of spaces}
\label{202002211550}

In this section, we address a characterization problem of the $d$-dimensional cospan categories of finite CW-spaces.
Recall some algebraic characterizations of the $d$-dimensional cospan categories of finite CW-spaces for $d = 0,1$ in subsection \ref{202201180923}.
Note that such a characterization becomes more difficult for higher $d$.
We modify the problem to a characterization problem {\it relative to} the homotopy category of finite CW-spaces.
It yields a stronger version of our main theorem.

Here is a heuristic argument.
For $d =\infty$, it is easy to show that the category $\mathsf{Cosp}^{\simeq}_{\leq (d+1)} ( \mathsf{CW}^\mathsf{fin}_\ast )$ is generated by the homotopy category of spaces and its opposite category : every morphism $[ \Lambda ]$ is decomposed into $\iota^- (f_1) \circ \iota^+ ( f_0 )$ where $\Lambda = ( K_0 \stackrel{f_0}{\to} L \stackrel{f_1}{\leftarrow} K_1 )$.
We carefully enhance this observation to a universality.
In the case $d < \infty$, such a decomposition is not valid since $\dim L$ could be $(d+1)$ essentially even up to a homotopy.
We weaken the universality to a unification property and show that the cospan categories fulfill the unification property relative to some functors.
It is the most nontrivial part in our discussion.
From now on, we use notations $\mathcal{H}_d = \mathsf{Ho} ( \mathsf{CW}^\mathsf{fin}_{\ast, \leq d} )$ for an extended natural number $d = 0 , 1, \cdots , \infty$ and $\mathcal{H} = \mathcal{H}_\infty$ for short.

\subsection{Mayer-Vietoris positive-and-negative diagram}
\label{202201172154}

In category theory, a functor $F : \mathcal{D} \to \mathcal{S}$ is called a $\mathcal{D}$-shaped diagram in $\mathcal{S}$.
In a similar fashion only here, we call a pair of functors $F^+ : \mathcal{D} \to \mathcal{S}$ and $F^- : \mathcal{D}^\mathrm{op} \to \mathcal{S}$ by a {\it $\mathcal{D}$-shaped positive-and-negative diagram in $\mathcal{S}$} : $F^+$ and $F^-$ play a role of a positive diagram and a negative diagram due to the orientation depicted in the $\mathcal{S}$.
In this paper, we are mainly interested in a specific positive-and-negative diagram as follows.

\begin{Defn}
For a category $\mathcal{S}$, a {\it Mayer-Vietoris Positive-and-Negative diagram in $\mathcal{S}$ for pointed finite CW-spaces with dimension lower than or equal to $d$} , or a {\it MV PN diagram in $\mathcal{S}$ with $\dim \leq d$} for short, is a pair of functors $F^+ : \mathcal{H}_d \to \mathcal{S}$ and $F^- : \mathcal{H}_d^\mathrm{op} \to \mathcal{S}$.
It assigns the same object in $\mathcal{S}$ to any $K \in \mathcal{H}_d$, i.e. $F^+ (K) = F^- (K)$.
For a homotopy pushout diagram (\ref{202201132212}) in $\mathcal{H}_d$, we have (\ref{202201301150}).
\begin{equation}
\label{202201132212}
\begin{tikzcd}[sep=scriptsize]
L_0 \ar[r, "g_1"] & S \\
K \ar[u, "f_0"] \ar[r, "g_0"'] & L_1 \ar[u, "f_1"]
\end{tikzcd}
\end{equation}
\begin{equation}
\label{202201301150}
F^- (f_1) \circ F^+ (g_1)  = F^+ ( g_0 ) \circ F^- ( f_0 ) .
\end{equation}
For simplicity, we use a notation $F^\pm$ or the following :
\begin{align}
\notag
F^+ : \mathcal{H}_d \longrightarrow \mathcal{S} \longleftarrow \mathcal{H}_d^\mathrm{op} : F^- .
\end{align}
Moreover, if $d = \infty$, then we call $F^\pm = (F^+ , F^-)$ a {\it MV PN diagram in $\mathcal{S}$} for short, i.e. we omit $\dim \leq \infty$.
It is depicted as follows.
\begin{align}
\notag
F^+ : \mathcal{H} \longrightarrow \mathcal{S} \longleftarrow \mathcal{H}^\mathrm{op} : F^- .
\end{align}
\end{Defn}

\begin{Example}
\label{202201141400}
For a morphism $f : K \to L$ in $\mathcal{H}_d$, let $\iota_d^+ (f)$ be the morphism in the category $\mathsf{Cosp}^\simeq_{\leq (d+1)} ( \mathsf{CW}^\mathsf{fin}_\ast )$ induced by the cospan $( K \stackrel{f}{\to} L \stackrel{\mathrm{id}}{\leftarrow} L)$.
The assignment $\iota_d^+$ induces a functor $\iota_d^+ : \mathcal{H}_d \to \mathsf{Cosp}^\simeq_{\leq (d+1)} ( \mathsf{CW}^\mathsf{fin}_\ast )$.
By using the dagger operation $\dagger$ on $\mathsf{Cosp}^\simeq_{\leq (d+1)} ( \mathsf{CW}^\mathsf{fin}_\ast )$, we define $\iota_d^- = \dagger \circ \iota_d^+ :  \mathcal{H}_d^\mathrm{op} \to \mathsf{Cosp}^\simeq_{\leq (d+1)} ( \mathsf{CW}^\mathsf{fin}_\ast )$.
We abbreviate $\iota_\infty^+ , \iota_\infty^-$ to $\iota^+, \iota^-$.
Then the following is a MV PN diagram in $\mathsf{Cosp}^\simeq_{\leq (d+1)} ( \mathsf{CW}^\mathsf{fin}_\ast )$ with $\dim \leq d$ by definitions :
\begin{align}
\notag
\iota_d^+ : \mathcal{H}_d \longrightarrow \mathsf{Cosp}^\simeq_{\leq (d+1)} ( \mathsf{CW}^\mathsf{fin}_\ast ) \longleftarrow \mathcal{H}_d^\mathrm{op} : \iota_d^- .
\end{align}
\end{Example}

\begin{Example}
\label{202201251638}
We give typical and instructive examples which arise from classical algebraic topology.
Consider $\mathcal{S} = \mathsf{Rel}_\ast$, i.e. the relation category of pointed sets.
The objects are pointed sets, and morphisms from $X$ to $Y$ are pointed subsets in $X \times Y$, i.e. pointed relations over $X$ and $Y$.
The identity is given by the diagonal sets and the composition is induced by the composition of relations.
Let $R^+$ and $R^-$ be the functors which assigns the graph relation and its transposition to a map respectively : for a pointed map $f : X \to Y$ we set
\begin{align}
R^+ (f) &= \{ (x , f(x)) \in X \times Y ~;~ x \in X \} , \notag \\
R^- (f) &= \{ ( f(x) , x ) \in Y \times X ~;~ x \in X \} . \notag
\end{align}
For a distinguished pointed space $X$, let $E : \mathcal{H}^\mathrm{op} \to \mathsf{Set}$ be a functor which assigns the homotopy set of maps $[ K , X ]$ to each $K$.
Let $F^+, F^-$ be compositions $F^+ = R^- \circ E^\mathrm{op}$ and $F^- = R^+ \circ E$.
Then $F^\pm$ gives a MV PN diagram in $\mathsf{Rel}_\ast$.
In fact, the equation $F^- (f_1) \circ F^+ (g_1)  = F^+ ( g_0 ) \circ F^- ( f_0 )$ in (\ref{202201301150}) is equivalent with the following proposition : for any $[a] \in [L_0 , X]$ and $[b] \in [L_1 , X]$,
\begin{align}
\notag
\exists [c] \in [S, X] \left( [c\circ g_1] = [a] , ~ [c \circ f_1] = [b] \right)
\Leftrightarrow
[a \circ f_0] = [b \circ g_0 ] .
\end{align}
Note that it is immediate from the homotopy pushout condition of (\ref{202201132212}).
\end{Example}

%%%%%%%%%%%%%%%%%%%%%%
\subsection{Universality of $\infty$-dimensional cospan category}
\label{202201301625}

In this subsection, we characterize the cospan categories of finite CW-spaces based on a universality defined below.

\begin{Defn}
A MV PN diagram $Z^\pm$in $\mathcal{U}$ (with $\dim \leq \infty$) satisfies a {\it universality} if the following condition holds :
\begin{itemize}
\item[]
For any MV PN diagram $F^\pm$ in a category $\mathcal{S}$, there exists a unique functor $F : \mathcal{U} \to \mathcal{S}$ such that the following diagram strictly commutes : 
\begin{equation}
\notag
\begin{tikzcd}
\mathcal{H} \ar[r, "Z^+"]  \ar[dr, "F^+"']& \mathcal{U} \ar[d, "F"] &  \mathcal{H}^\mathrm{op}  \ar[l, "Z^-"'] \ar[dl, "F^-"] \\
& \mathcal{S} &
\end{tikzcd}
\end{equation}
\end{itemize}
If such $F$ exists, then it is called a {\it unification of $F^+$ and $F^-$}.
Note that such a MV PN diagram satisfying the universality, if exists, uniquely exists up to a category isomorphism.
\end{Defn}

Recall $\iota^\pm = \iota^\pm_\infty$ in Example \ref{202201141400}.

\begin{theorem}
\label{202201172148}
The MV PN diagram $\iota^\pm$ satisfies the universality.
\end{theorem}

Before we prove the theorem, we introduce an important notion which is used in the remaining part of this paper.
Here we consider general $d$ not only $\infty$, and also spans of spaces not only cospans for later use.

\begin{Defn}
Let $F^\pm$ be a MV PN diagram in $\mathcal{S}$ with $\dim 
\leq d$.
For a $d$-dimensional cospan $\Lambda = ( K_0 \stackrel{f_0}{\to} L \stackrel{f_1}{\leftarrow} K_1 )$, we define a morphism $\int_{\Lambda} F^\pm$ from $F^+ (K_0) = F^- (K_0)$ to $F^+ (K_1) = F^- (K_1)$ in $\mathcal{S}$ by
\begin{align}
\notag
\int_\Lambda F^\pm \stackrel{\mathrm{def.}}{=}  F^- (f_1) \circ F^+ (f_0) .
\end{align}
Here, $\circ$ is the composition in $\mathcal{S}$.
Similarly for a span $\mathrm{V} = ( K_0 \stackrel{f_0}{\leftarrow} L \stackrel{f_1}{\to} K_1 )$ of spaces, we define a morphsim $\int_\mathrm{V} F^\pm \stackrel{\mathrm{def.}}{=}  F^+ (f_1) \circ F^- (f_0)$.
\end{Defn}
\begin{remark}
The notation $\int$ is motivated by the finite path-integral in DWFQ model and TVBW model where $\mathcal{S}$ is replaced by a relation category (or a category of Hilbert spaces) and $F^+, F^-$ are replaced by the finite path-integral and its adjoint (or dagger) respectively.
In spite of such motivation, our strategy is more general in the sense that we also deal with spans of spaces in subsubsection \ref{202201201217}.
\end{remark}

\begin{proof}[Proof of Theorem \ref{202201172148}]
Let $F^\pm$ be a MV PN diagram in $\mathcal{S}$.
Here we give a concrete construction of a unification $\hat{F}$ of $F^+$ and $F^-$ relative to $i_d$ as follows.
The functor $\hat{F}$ assigns $F^+ (K) = F^- (K)$ to every object $K$ of $\mathsf{Cosp}^{\simeq}(\mathsf{CW}^\mathsf{fin}_{\ast} )$.
Let $\Lambda = ( K_0 \stackrel{f_0}{\to} L \stackrel{f_1}{\leftarrow} K_1 )$ be a cospan in the homotopy category.
The functor $\hat{F}$ assigns $\int_\Lambda F^\pm$ to the morphism $[\Lambda ]$ in $\mathsf{Cosp}^{\simeq} (\mathsf{CW}^\mathsf{fin}_{\ast} )$.
It gives a well-defined functor by definitions.
Then one can directly prove that $\hat{F}$ is a unification of $F^+$ and $F^-$.

We prove the uniqueness of the unification.
Let $F$ be a unification of $F^\pm$.
Let $\Lambda = ( K_0 \stackrel{f_0}{\to} L \stackrel{f_1}{\leftarrow} K_1 )$ be a cospan of pointed finite CW-spaces.
Then we have $F ( [ \Lambda ] ) = F ( \iota^- ( f_1 ) ) \circ F ( \iota^+ ( f_0 ) ) = F^- ( f_1 ) \circ F^+ ( f_0 ) = \hat{F} ( [ \Lambda ] )$.
Thus, a unification $F$ of $F^+$ and $F^-$ is uniquely determined by $F^\pm$.
It completes the proof.
\end{proof}

\begin{Example}
The unification of $\iota^+$ and $\iota^-$ in Example \ref{202201141400} is the identity functor.
\end{Example}

\begin{Example}
We obtain a unification $F : \mathsf{Cosp}^{\simeq}(\mathsf{CW}^\mathsf{fin}_{\ast} ) \to \mathsf{Rel}_\ast$ which is characterized by $F \circ \iota^+ = R^- \circ E^\mathrm{op}$ and $F \circ \iota^- = R^+ \circ E$ by applying the theorem to Example \ref{202201251638}.
\end{Example}

%%%%%%%%%%%%%%%%%%%%%%
\subsection{Unification property of finite-dimensional cospan categories}
\label{202201301626}

In this subsection, we deal with a characterization problem of $ \mathsf{Cosp}^\simeq_{\leq d} ( \mathsf{CW}^\mathsf{fin}_\ast )$ (cf. Theorem \ref{202201172148}) when $d$ is finite.
It is not obvious whether it is possible to characterize the category by a universality, but it turns out that it fulfills a weak universality relative to some canonical functors.
We call such a weak universality by a unification property as follows :

\begin{Defn}
Let $d = 0, 1, \cdots , \infty$ be an extended natural number.
Let $P : \mathcal{H}_d \to \mathcal{H}_{d+1}$ be a functor.
The MV PN diagram $Z^\pm$ in a category $\mathcal{U}$ satisfies a {\it unification property relative to $P$} if the following condition holds :
\begin{itemize}
\item[]
For any MV PN diagram $F^\pm$ in a category $\mathcal{S}$ with $\dim \leq (d+1)$, there exists a unique functor $F : \mathcal{U} \to \mathcal{S}$ such that the following diagram strictly commutes : 
\begin{equation}
\notag
\begin{tikzcd}
\mathcal{H}_d \ar[r, "\iota_d^+"] \ar[d, "P"'] & \mathcal{U} \ar[d, "F"] &  \mathcal{H}_d^\mathrm{op}  \ar[l, "\iota_d^-"']  \ar[d, "P^\mathrm{op}"] \\
\mathcal{H}_{d+1}  \ar[r, "F^+"']  & \mathcal{S} & \mathcal{H}_{d+1}^\mathrm{op} \ar[l, "F^-"]
\end{tikzcd}
\end{equation}
If such $F$ exists, then it is called a {\it unification of $F^+$ and $F^-$ relative to $P$} and denoted by $P^\cup ( F^\pm )$.
\end{itemize}
\end{Defn}

\begin{remark}
The prefix {\it weak} means that it does not characterize the category.
However, such a unification property yields a canonical construction of TQFT's \cite{kim2021pair}.
\end{remark}

Recall the suspension of pointed spaces.
For a pointed space $K$, we set the suspension by the smash product with the pointed circle $S^1_\ast$,
\begin{align}
\notag
\Sigma K = S^1_\ast \wedge K .
\end{align}
If $K$ is a finite CW-space, then so is $\Sigma K$.
Moreover, if $\dim K \leq d$, then $\dim K \leq (d+1)$.
Denote by $\Sigma_d : \mathcal{H}_d \to \mathcal{H}_{d+1}$ the assignment of suspensions.

\begin{theorem}
\label{202201172145}
Let $d = 0, 1, \cdots , \infty$ be an extended natural number.
The MV PN diagram $\iota_d^\pm$ satisfies the unification property relative to
\begin{itemize}
\item
the embedding functor $i_d : \mathcal{H}_d \to \mathcal{H}_{d+1}$,
\item
and the suspension functor $\Sigma_d : \mathcal{H}_d \to \mathcal{H}_{d+1}$.
\end{itemize}
Furthermore, the unification and the pullback commutes in the following sense :
Let $F^\pm$ be a MV-PN diagram in a category $\mathsf{S}$ with $\dim \leq (d+1)$, and $j_d : \mathsf{Cosp}^\simeq_{\leq d} ( \mathsf{CW}^\mathsf{fin}_\ast ) \to \mathsf{Cosp}^\simeq_{\leq (d+1)} ( \mathsf{CW}^\mathsf{fin}_\ast )$ be the embedding functor.
If $d \geq 1$, then the unifications of $F^\pm$ relative to $\Sigma_d$ and of $\Sigma_d^\ast ( F^\pm )$ relative to $i_{d-1}$ exist, and we have
\begin{align}
\notag
j_d^\ast ( \Sigma_d^\cup ( F^\pm ) ) = i_{d-1}^\cup ( \Sigma_d^\ast ( F^\pm ) ) .
\end{align}

\end{theorem}

\begin{remark}
The unification property relative to the embedding functor for $d = \infty$ obviously implies Theorem \ref{202201172148}.
\end{remark}

\subsubsection{The unification property relative to the embedding functor}

\begin{Lemma}
\label{201912231215}
Let $d = 0, 1, \cdots$ be a natural number.
The category $\mathsf{Cosp}^{\simeq}_{\leq (d+1)} (\mathsf{CW}^\mathsf{fin}_{\ast} )$ is generated by cospans $\dim_\mathrm{c} \Lambda \leq d$.
\end{Lemma}
\begin{proof}
We give a proof of a stronger statement :
For a cospan $\Lambda$ with $\dim \Lambda \leq (d+1)$, there exist composable cospans $\Lambda_0, \Lambda_1, \Lambda_2$ of pointed finite CW-complexes subject to following conditions.
\begin{enumerate}
\item
$\dim_\mathrm{c} \Lambda_j \leq d$.
\item
The cospans $\Lambda_2, \Lambda_1, \Lambda_0$ are composable, i.e. the target of $\Lambda_0$ ($\Lambda_1$, resp.) and the source of $\Lambda_1$ ($\Lambda_2$, resp.) are the same.
\item
We have $[ \Lambda ] = [ \Lambda_2 ] \circ [ \Lambda_1 ] \circ [ \Lambda_0 ]$ in the category $\mathsf{Cosp}^{\simeq}_{\leq (d+1)} (\mathsf{CW}^\mathsf{fin}_{\ast} )$.
\end{enumerate}
First, we show the claim for a special case.
We assume that $\Lambda = ( K \stackrel{f}{\to} L \leftarrow \mathrm{pt} )$ where $\dim K \leq d$ and $\dim L \leq (d+1)$.
Recall that $f$ in the diagram $( K \stackrel{f}{\to} L \leftarrow \mathrm{pt} )$, in fact, represents the homotopy equivalence class of $f$.
We choose pointed finite CW-complex structures $X_L$ on $L$ and $X_K$ on $K$.
We also choose a representative $f$ which is a cellular map with respect to $X_L, X_K$.
Then we have $f(K) \subset L^{(d)}$ where $L^{(d)}$ is the $d$-skeleton of $X_L$.
Denote by $f^\prime : K \to L^{(d)}$ the induced map.
Let $\varphi_j : D^{d+1} \to L,~~r = 1,2, \cdots, k$ be characteristic maps of $(d+1)$-cells of $X_L$.
Let $\psi : \vee_j ( S^{d})^+  \to L^{(d)}$ be the pointed map induced by the wedge sum $\vee_j (\varphi_j |_{S^{d}})$.
Let $c : \vee_j (S^{d})^+ \to  \vee_j \mathrm{pt}^+$ be the pointed map induced by the collapsing map from $S^{d}$ to $\mathrm{pt}$.
Then $L = L^{(d+1)}$ is homeomorphic to the double mapping cylinder of $\psi$ and $c$.
Hence, we have $[ \Lambda_1] \circ [ \Lambda_0] = [ \Lambda ]$ where $\Lambda_0 = \left( K \stackrel{f^\prime}{\to} L^{(d)} \stackrel{\psi}{\leftarrow} \vee_j ( S^{d})^+ \right)$ and $\Lambda_1 = \left( \vee_j ( S^{d})^+ \stackrel{c}{\to} \vee_j  \mathrm{pt}^+ \leftarrow \mathrm{pt} \right)$.
Furthermore, the components of $\Lambda_0$ and $\Lambda_1$ have dimension lower than or equal to $d$ since $d \geq 0$.

Now we give a general proof.
Any object $K$ in $\mathsf{Cosp}^{\simeq}_{\leq (d+1)} (\mathsf{CW}^\mathsf{fin}_{\ast} )$, i.e. a finite CW-space with dimension $\leq d$, has a self-dual given by the canonical cospans $( K \vee K \to K \leftarrow \mathrm{pt} )$ and $( \mathrm{pt} \to K \leftarrow K \vee K )$.
Denote by $\mathrm{ev}_K$ and $\mathrm{coev}_K$ the induced morphisms in $\mathsf{Cosp}^{\simeq}_{\leq (d+1)} (\mathsf{CW}^\mathsf{fin}_{\ast} )$.
Then we have $[\Lambda] = ( \mathrm{ev}_{K_1} \vee \mathrm{id}_{K_1} ) \circ ( [\Lambda ] \vee \mathrm{id}_{K_1} \vee \mathrm{id}_{K_1}) \circ ( \mathrm{id}_{K_0} \vee \mathrm{coev}_{K_1})$.
See Figure \ref{201912182123}.
Note that $( \mathrm{ev}_{K_1} \vee \mathrm{id}_{K_1} ) \circ ( [\Lambda ] \vee \mathrm{id}_{K_1} \vee \mathrm{id}_{K_1}) = ( \mathrm{ev}_{K_1} \circ ([\Lambda ] \vee \mathrm{id}_{K_1}) ) \vee \mathrm{id}_{K_1}$.
By the previous discussion, there exist cospans $\Lambda_2, \Lambda_1$ such that $( \mathrm{ev}_{K_1} \circ ([\Lambda ] \vee \mathrm{id}_{K_1} ) ) = [ \Lambda_2] \circ [\Lambda_1 ]$ and  $\dim_\mathrm{c} \Lambda_2, \dim_\mathrm{c} \Lambda_1 \leq d$.
Moreover, $\Lambda_0 = ( \mathrm{id}_{K_0}  \vee \mathrm{coev}_{K_1})$ satisfies $\dim_\mathrm{c} \Lambda_0 \leq d$.
By the choices, we have $[\Lambda] = [\Lambda_2] \circ [\Lambda_1] \circ [\Lambda_0]$.
It completes the proof.
\begin{figure}[h!]
  \includegraphics[width=8.6cm]{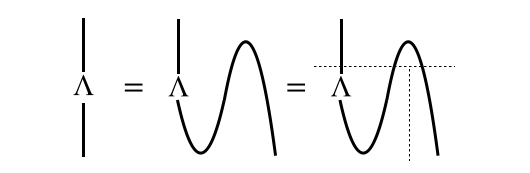}
  \caption{A decomposition of a cospan $\Lambda$}
  \label{201912182123}
\end{figure}
\end{proof}

\begin{proof}[Proof of the first part of Theorem \ref{202201172145}]
Let $F^\pm$ be a MV-PN diagram in a category $\mathsf{S}$ with $\dim \leq (d+1)$.
We apply the construction of $\hat{F}$ in Theorem \ref{202201172148} with careful consideration of dimensions to obtain a unification.
It suffices to prove the uniqueness of the unification.
Let $F$ be a unification of $F^\pm$.
Let $\Lambda = ( K_0 \stackrel{f_0}{\to} L \stackrel{f_1}{\leftarrow} K_1 )$ be a cospan of pointed finite CW-spaces.
Then we have $F ( [ \Lambda ] ) = F ( \iota^- ( f_1 ) ) \circ F ( \iota^+ ( f_0 ) ) = F^- ( f_1 ) \circ F^+ ( f_0 ) = \hat{F} ( [ \Lambda ] )$.
Hence, $F ( [ \Lambda ] )$ is  determined by the given $F^\pm$ if $\Lambda$ is a cospan whose components have dimension $\leq d$.
By Lemma \ref{201912231215}, all the morphisms in $\mathsf{Cosp}^{\simeq}_{\leq (d+1)} (\mathsf{CW}^\mathsf{fin}_{\ast} )$ is decomposed into some morphisms of the above type.
Thus, a unification $F$ of $F^+$ and $F^-$ relative to $i_d$ is uniquely determined by $F^\pm$.
\end{proof}

\subsubsection{The unification property relative to the suspension functor}
\label{202201201217}

\begin{Defn}
\label{201912241125}
Denote by $\bar{t} \in S^1$ the complex conjugate of $t \in S^1 \subset \mathbb{C}$.
We define a conjugate $\tau_K : \Sigma K \to \Sigma K ; [t,x] \mapsto [\bar{t} , x]$ where $[t, x]$ is the induced element of $\Sigma K$ for elements $t \in S^1$ and $x \in K$.
Let $\Lambda = \left( K_0 \stackrel{f_0}{\to} L \stackrel{f_1}{\leftarrow} K_1 \right)$ be a cospan in $\mathsf{Top}_\ast$.
Denote by $p_0 , p_1$ the collapsing maps from the mapping cone $\mathrm{Cone} ( f_0 \vee f_1 )$ to the suspensions $\Sigma K_0$ and $\Sigma K_1$ respectively.
We define a span $\mathrm{T} \Sigma ( \Lambda)$ in $\mathsf{Top}_\ast$ by
\begin{align}
\notag
\mathrm{T} \Sigma ( \Lambda)
\stackrel{\mathrm{def.}}{=} \left( \Sigma K_0 \stackrel{\tau_{K_0} \circ p_0}{\longleftarrow} \mathrm{Cone} ( f_0 \vee f_1 ) \stackrel{p_1}{\longrightarrow} \Sigma K_1 \right) .
\end{align}
\end{Defn}

\begin{remark}
Note that the mapping cone $\mathrm{Cone} ( f_0 \vee f_1 )$ might not be a CW-space even if so are $K_0,K_1,L$.
A cospan $\Lambda$ in $\mathsf{CW}_\ast$ is {\it cellular} if $f_0,f_1$ are cellular maps for some CW-complex structures on $K_0,K_1,L$.
Note that for any cospan $\Lambda$ in $\mathsf{CW}_\ast$ there exists cellular cospan $\Lambda^\prime$ in $\mathsf{CW}_\ast$ such that $\Lambda \simeq \Lambda^\prime$ by the cellular approximation theorem.
We shall modify any cospans to such cellular ones up to a homotopy equivalence.
It is easy to show that $\mathrm{T} \Sigma ( \Lambda)$ is a span in $\mathsf{CW}_\ast$ if a cospan $\Lambda$ in $\mathsf{CW}_\ast$ is cellular.
Furthermore, if $\dim \Lambda \leq d$, then $\dim_\mathrm{c} \mathrm{T} \Sigma ( \Lambda) \leq d$.
\end{remark}

\begin{proof}[Proof of the second part of Theorem \ref{202201172145}]
We decompose a construction of a unification $\Sigma_d^\cup ( F^\pm ) = \check{F}$ to lemmas.
The functor $\check{F}$ assigns $F^+ ( \Sigma K) = F^- ( \Sigma K )$ to an object $K$ of $\mathsf{Cosp}^{\simeq}_{\leq d} (\mathsf{CW}^\mathsf{fin}_{\ast} )$.
It is well-defined since the domain of $F^+,F^-$ consists of pointed finite CW-spaces $X$ with $\dim X \leq d$.
The functor $\check{F}$ assigns $\int_{\mathrm{T}\Sigma (\Lambda)} F^\pm = F^+ ( p_1 ) \circ F^- ( \tau_{K_0} \circ p_0 )$ to a morphism $[\Lambda]$ where $\Lambda = ( K_0 \stackrel{f_0}{\to} L \stackrel{f_1}{\leftarrow} K_1 )$.
To verify that it is well-defined, we show that $\int_{\mathrm{T}\Sigma (\Lambda)} F^\pm$ depends on the homotopy equivalence class of $\Lambda$.
See Lemma \ref{202111031525}.
Moreover, we prove that the assignment preserves the composition by Lemma \ref{201912231233}.
All that remain is to show that the functor $\check{F}$ is a unification of $F^+$ and $F^-$ relative to $\Sigma_d$.
It suffices to prove that $\check{F} ( \iota^+ ( f )) = F^+ (\Sigma_d f)$ for a map $f: K \to L$ lying in $\mathcal{H}_d$.
Denote by $p_0,p_1$ the collapsing maps from the mapping cone $\mathrm{Cone} ( f \vee \mathrm{id}_L )$ to $\Sigma K$ and $\Sigma L$ respectively.
Note that the map $p_0$ is a pointed homotopy equivalence, especially so $\tau_K \circ p_0$ is.
Then $\check{F} ( \iota^+ ( f )) = F^+ (\Sigma_d f)$ follows from the diagram below.
The uniqueness of the unification is proved similarly with the proof of the first part.
\begin{equation}
\notag
\begin{tikzcd}[sep=scriptsize]
& \mathrm{Cone} (f \vee \mathrm{id}_L ) \ar[dl, "\tau_K \circ p_0"'] \ar[dr, "p_1"] \ar[dd, "\tau_K \circ p_0"] & \\
\Sigma K  & & \Sigma L  \\
& \Sigma K \ar[ur, "\Sigma f"'] \ar[ul, "\mathrm{id}_{\Sigma K}"] & 
\end{tikzcd}
\end{equation}
\end{proof}

From now on, we prove the lemmas appearing before.

\begin{Lemma}
\label{202111031525}
Let $\Lambda, \Lambda^\prime$ be cospans in $\mathsf{Top}_\ast$.
A homotopy equivalence of cospans $\Lambda \simeq \Lambda^\prime$ implies a homotopy equivalence of spans $\mathrm{T} \Sigma ( \Lambda) \simeq \mathrm{T} \Sigma ( \Lambda^\prime )$.
\end{Lemma}
\begin{proof}
For a map $f : K \to L$, let $p_f : \mathrm{Cone} ( f ) \to \Sigma K$ and $j_f : L \to \mathrm{Cone} ( f )$ be the canonical maps associated with $f$.
Consider a map $h : L \to L^\prime$. 
Let $g_{f,h} : \mathrm{Cone} ( f ) \to \mathrm{Cone} ( h \circ f)$ be the pointed map which is induced by gluing the identity on the cone of $K$ and $h : L  \to L^\prime$.
We have an obvious homotopy pushout diagram as follows :
\begin{equation}
\notag
\begin{tikzcd}[sep=scriptsize]
\mathrm{Cone} (f) \ar[r, "g_{f,h}"] & \mathrm{Cone} ( h \circ f) \\
L \ar[r, "h"'] \ar[u, "j_f"] & L^\prime \ar[u, "j_{h\circ f}"']
\end{tikzcd}
\end{equation}
Especially, if $h$ is a homotopy equivalence, then so is $g_{f,h}$.
Moreover, we have $p_{h\circ f} \circ g_{f,h} = p_f$.
Now suppose that $\Lambda = ( K_0 \stackrel{f_0}{\to} L \stackrel{f_1}{\leftarrow} K_1 )$, $\Lambda^\prime = ( K_0 \stackrel{g_0}{\to} L^\prime \stackrel{g_1}{\leftarrow} K_1 )$ and a homotopy equivalence $h : L \to L^\prime$ induces $\Lambda \simeq \Lambda^\prime$.
Let $\Lambda^{\prime\prime} = ( K_0 \stackrel{h \circ f_0}{\to} L^\prime \stackrel{h \circ f_1}{\leftarrow} K_1 )$.
Then for $K = K_0 \vee K_1$ and $f = f_0 \vee f_1 : K_0 \vee K_1 \to L$, the homotopy equivalence $g_{f,h}$ gives $\mathrm{T} \Sigma ( \Lambda) \simeq \mathrm{T} \Sigma ( \Lambda^{\prime\prime} )$.
Moreover, $\mathrm{T} \Sigma ( \Lambda^\prime ) \simeq \mathrm{T} \Sigma ( \Lambda^{\prime\prime} )$ follows from the elementary fact that $u_0,u_1 : K \to L$ are homotopic with each other, then there exists a homotopy equivalence $k : \mathrm{Cone} (u_0) \to \mathrm{Cone} (u_1)$ such that $p_{u_1} \circ k \simeq p_{u_0}$.
We apply the fact to $u_0 = h \circ (f_0 \vee f_1)$ and $u_1 = g_0 \vee g_1$.
It completes the proof that $\mathrm{T} \Sigma ( \Lambda) \simeq \mathrm{T} \Sigma ( \Lambda^\prime )$.
\end{proof}

\begin{Lemma}
\label{201912231233}
Let $\Lambda, \Lambda^\prime, \Lambda^{\prime\prime}$ be $(d+1)$-dimensional cospans of pointed finite CW-spaces such that $[\Lambda^{\prime\prime} ] = [ \Lambda^\prime ] \circ [ \Lambda ]$.
For a MV PN diagram $F^\pm$ with $\dim \leq (d+1)$, we have
\begin{align}
\notag
\left(\int_{\mathrm{T}\Sigma (\Lambda^\prime)} F^\pm \right) \circ  \left( \int_{\mathrm{T} \Sigma (\Lambda)} F^\pm \right) =  \int_{\mathrm{T}\Sigma ( \Lambda^{\prime\prime} )} F^\pm .
\end{align}
\end{Lemma}
\begin{proof}
Let $\Lambda = ( K_0 \stackrel{f_0}{\to} L \stackrel{f_1}{\leftarrow} K_1 )$ and $\Lambda^\prime = ( K_1 \stackrel{f^\prime_1}{\to} L^\prime \stackrel{f^\prime_2}{\leftarrow} K_2 )$.
Without loss of generality, one could choose a representative $\Lambda^{\prime\prime} = ( K_0 \stackrel{f^{\prime\prime}_0}{\to} L^{\prime\prime} \stackrel{f^{\prime\prime}_2}{\leftarrow} K_2 )$ as the homotopy pushout composition.
Especially, $L^{\prime\prime} = \mathrm{Cyl} (f_1 , f^\prime_1)$.
Denote by $p_0,p_1$ the collapsing maps from the mapping cone $\mathrm{Cone} (f_0 \vee f_1)$ to the suspensions $\Sigma K_0,\Sigma K_1$, and by $p^\prime_1, p^\prime_2$ the collapsing maps from the mapping cone $\mathrm{Cone} (f^\prime_1 \vee f^\prime_2)$ to the suspensions $\Sigma K_1,\Sigma K_2$.
Denote by $q_0, q_2$ the collapsing maps from the mapping cone $\mathrm{Cone} (f^{\prime\prime}_0 \vee f^{\prime\prime}_2 )$ to the mapping cones $\mathrm{Cone} (f_0 \vee f_1)$ and $\mathrm{Cone} (f^\prime_1 \vee f^\prime_2)$ respectively.
It turns out that the diagram (\ref{201912232217}) commutes up to a homotopy.
Recall that we use the conjugate $\tau_{K_0}$ in Definition \ref{201912241125} so that this diagram to commute.
Moreover the square diagram (\ref{201912232217}) is a homotopy pushout diagram.
\begin{equation}
\label{201912232217}
\begin{tikzcd}[sep=scriptsize]
\mathrm{Cone} (f_0 \vee f_1)  \ar[r, "p_1"] &  \Sigma K_1  \\
\mathrm{Cone} (f^{\prime\prime}_0 \vee f^{\prime\prime}_2 ) \ar[r, "q_2"] \ar[u, "q_0"] & \mathrm{Cone} (f^\prime_1 \vee f^\prime_2) \ar[u, "\tau_{K_1} \circ p^\prime_1"']
\end{tikzcd}
\end{equation}
It gives $F^-(\tau_{K_1} \circ p^\prime_1) \circ F^+(p_1) = F^+(q_2) \circ F^-(q_0)$.
It completes the proof.
\end{proof}

\subsubsection{Commutativity of the unification and the pullback}

Finally, we prove the last part of Theorem \ref{202201172145}.
It follows from the following lemma.
As notations, we set $\Sigma ( \Lambda ) = ( \Sigma K_0 \stackrel{\Sigma f_0}{\to} \Sigma L \stackrel{\Sigma f_1}{\leftarrow} \Sigma K_1 )$ for $\Lambda = ( K_0 \stackrel{f_0}{\to} L \stackrel{f_1}{\leftarrow} K_1 )$.

\begin{Lemma}
\label{202112111313}
Let $F^\pm$ be a MV PN diagram in $\mathcal{S}$ with $\dim \leq (d+1)$.
For a cospan $\Lambda$ with $\dim_\mathrm{c} \Lambda \leq d$, we have $\dim_\mathrm{c} \Sigma ( \Lambda ) \leq (d+1)$ and 
\begin{align}
\notag
\int_{\mathrm{T}\Sigma ( \Lambda )} F^\pm  = \int_{\Sigma ( \Lambda )} F^\pm  .
\end{align}
\end{Lemma}
\begin{proof}
The statement is immediate from the fact that the diagram below is a homotopy pushout diagram.
\begin{equation}
\notag
\begin{tikzcd}[sep=scriptsize]
\Sigma K_0 \ar[r, "\Sigma f_0"] & \Sigma L \\
\mathrm{Cone} ( f_0 \vee f_1 ) \ar[r , "p_1"'] \ar[u, "\tau_{K_0} \circ p_0"] &  \Sigma K_1 \ar[u, "\Sigma f_1"']
\end{tikzcd}
\end{equation}
\end{proof}

\begin{remark}
Let $\dim K_0 , \dim K_1 \leq d$.
Note that the right hand side in Lemma \ref{202112111313} could be defined only if $\dim L \leq d$ whereas the left hand side is defined even if $\dim L = d$.
\end{remark}

%%%%%%%%%%%%%%%%%%%%%%%%%%%%%%%%%%%%%%%%%

\section{A pair of extensions of Mayer-Vietoris functors}

In this section, we first give a construction of MV PN diagrams starting from some functors whose codmain is an abelian category $\mathcal{A}$.
Those functors are subject to the Mayer-Vietoris condition and named {\it Mayer-Vietoris functors}.
We obtain MV PN diagrams whose target is a relation category associated with $\mathcal{A}$, which was classically studied.
We give an overview of those relation categories in subsection \ref{202201301328}.
The goal of this section is to prove our main theorem by applying the unification property in Theorem \ref{202201172145}.

\subsection{Mayer-Vietoris functors}
\label{202201301628}

In this subsection, we define $\mathcal{A}$-valued Mayer-Vietoris functors and give some examples.

\begin{Defn}
\label{202202021632}
For an extended natural number $d \in \mathbb{N} \cup \{ \infty \}$, an {\it $\mathcal{A}$-valued Mayer-Vietoris functor for finite CW-spaces with dimension $\leq d$}, or an {\it $\mathcal{A}$-valued MV functor with $\dim \leq d$} for short, is a functor $E : \mathcal{H}_d \to \mathcal{A}$ subject to the following condition.
Consider a homotopy pushout diagram in $\mathcal{H}_d$, say (\ref{202201132212}).
The induced chain complex (\ref{202201172142}) is exact.
If $d= \infty$, then we abbreviate the terminology to a $\mathcal{A}$-valued MV functor.
\begin{align}
\label{202201172142}
E(K) \stackrel{( (f_0)_\ast , -(g_0)_\ast)}{\longrightarrow} E(L_0) \oplus E(L_1) \stackrel{(g_1)_\ast + (f_1)_\ast}{\longrightarrow} E(S) .
\end{align}
An $\mathcal{A}$-valued MV functor with $\dim \leq d$ is {\it pointed} if the object $E( \ast )$ corresponding to the one-point space is a zero object of $\mathcal{A}$.
\end{Defn}

\begin{remark}
\label{202205261256}
The definition of $\mathcal{A}$-valued Brown functors in \cite{kim2021pair} corresponds to that of {\it pointed} $\mathcal{A}$-valued Mayer-Vietoris functors here.
We reintroduce the terminology to avoid a confusion caused by the classical meaning of Brown functor : a Brown functor is given by a functor from the opposite homotopy category of pointed {\it connected} CW-complexes.
\end{remark}

\begin{remark}
\label{202201301359}
Any pointed CW-spaces $K,L$ induce a homotopy pushout diagram as follows.
\begin{equation}
\notag
\begin{tikzcd}[sep=scriptsize]
K \vee L \ar[r] & L \\
K \ar[r] \ar[u] & \ast \ar[u]
\end{tikzcd}
\end{equation}
It turns out that the associated exact sequence is split so that we have a natural isomorphism $E(\ast) \oplus E(K \vee L) \cong E(K) \oplus E(L)$.
In particular, an $\mathcal{A}$-valued MV functor is pointed if and only if $E : \mathcal{H}_d \to \mathcal{A}$ is enhanced to a symmetric monoidal functor.
Here the symmetric monoidal structure on $\mathcal{H}_d$ and $\mathcal{A}$ is induced by the wedge sum of spaces and the biproduct respectively.
\end{remark}

\begin{remark}
In \cite{kim2021pair}, we defined a terminology, {\it $\mathcal{A}$-valued $d$-dimensional Brown functor}.
It is a functor $E : \mathcal{H}_d \to \mathcal{A}$ which maps a wedge sum to a biproduct in $\mathcal{A}$, and a square diagram approximated by a triad of spaces to an exact square diagram in $\mathcal{A}$.
By Remark \ref{202201301359}, it coincides with an $\mathcal{A}$-valued MV functor with $\dim \leq d$ which is pointed.
In this sense, we reintroduce the terminology in order to avoid a confusion caused by the classical meaning of Brown functor.
\end{remark}

\begin{Example}
Let $E$ be an $\mathcal{A}$-valued MV functor with $\dim \leq d$.
For an arbtirary object $A$ of $\mathcal{A}$, a functor induced by $E^\prime (K) = E (K) \oplus A$ gives an $\mathcal{A}$-valued MV functor with $\dim \leq d$.
\end{Example}

\begin{Example}
Let $\widetilde{E}_\bullet$ be a generalized homology theory, i.e. a sequence of functors $\widetilde{E}_q : \mathbb{CW}_\ast \to \mathsf{Ab}, q \in \mathbb{Z}$ which is equipped with the suspension isomorphisms $\widetilde{E}_{q+1} \circ \Sigma \cong \widetilde{E}_q$ and satisfy the Eilenberg-Steenrod axiom except the dimension axiom.
Then the assignment of $\widetilde{E}_q (K)$ to $K$ gives a pointed $\mathsf{Ab}$-valued MV functor (with $\dim \leq \infty$).
Here $\mathsf{Ab}$ could be replaced by the category of $R$-modules $\mathsf{Mod}_R$.
There are similar examples arising from generalized cohomology theories.
\end{Example}

\begin{Example}
There is a formal generalization of generalized (co)homology theories to those valued at an arbitrary abelian category $\mathcal{A}$.
The only difference is that the functor $\widetilde{E}_q$ is valued at $\mathcal{A}$, i.e. the corresponding $\widetilde{E}_q (K)$ is an object of $\mathcal{A}$ and all the structure homomorphisms lie in $\mathcal{A}$.
Then $\mathcal{A}$-valued generalized (co)homology theories induce pointed $\mathcal{A}$-valued MV functors as above.
\end{Example}

\begin{Example}
We give an example such that $d$ is not possibly $\infty$.
Suppose that $\mathcal{A}$ is equipped with a volume \cite{kim2021integrals}.
For an abelian monoid $M$, an $M$-valued volume $v$ on $\mathcal{A}$ is an additive invariant of objects in $\mathcal{A}$ which is compatible with exact sequences : if there is an exact sequence $0 \to A \to B \to C \to 0$, then we have $v(B) = v(A) \cdot v(C)$.
For example, the order is a $M_1$-valued volume on the category of abelian groups where $M_1 = \mathbb{Q}^\times \cup \{ \infty \}$ is the multiplicative abelian monoid ; the dimension is $M_2$-valued volume on the category of vector spaces where $M_1 = \mathbb{Z} \cup \{ \infty \}$ is the additive abelian monoid.
There is also an example obtained from bicommutative Hopf algebras \cite{kim2020homology}.
Let $\mathcal{A}^\mathsf{fin}$ be the subcategory consisting of objects with a finite volume, i.e. an object $A$ such that $v(A)$ is invertible in $M$.
Let $\widetilde{E}_\bullet$ be a generalized homology theory valued at $\mathcal{A}$.
If $\widetilde{E}_q (S^0_\ast) \in \mathcal{A}^\mathsf{fin}$ where $S^0_\ast$ is a pointed $0$-sphere, then an extended natural number $d = d ( \widetilde{E}_\bullet ; q )$ is defined (Definition 4.26 \cite{kim2021pair}) so that $\widetilde{E}_q (K) \in \mathcal{A}^\mathsf{fin}$ for any $K \in \mathcal{H}_d$.
Here $d$ is possibly finite.
Then $\widetilde{E}_q$ induces a pointed $\mathcal{A}^\mathsf{fin}$-valued MV functor with $\dim \leq d$.
The readers are carefully referred to subsection 4.1 \cite{kim2021pair} for more examples and explanation.
\end{Example}

\subsection{Relation category of abelian categories}
\label{202201301328}

The relation category in Example \ref{202201251638} is defined by using relations over given two sets.
There is an additive generalization.
For abelian groups $X,Y$, an (additive) relation over $X$ and $Y$ is a subgroup of $X \times Y$.
In a parallel way, one can define a relation category of abelian groups $\mathsf{Rel} ( \mathsf{Ab} )$.
Based on the idea that a subgroup of $X$ could be replaced by an equivalence class of a monomorphism to $X$, i.e. a subobject of $X$, relations in an abelian category had been systematically studied \cite{mac1961algebra} \cite{meisen1974relations} \cite{puppe1962korrespondenzen} \cite{hilton1966correspondences} ; see \cite{gran2021introduction} for a recent introduction in a more general context.
In this subsection, we give an overview based on Hilton.
Fix an abelian category $\mathcal{A}$.
We assume that for any object $A$ all the subobjects of $A$ forms a set.
See Remark \ref{202202031403} for this hypothesis.

\begin{Defn}
Let $\mathrm{V}_0$, $\mathrm{V}_1$ be spans in $\mathcal{A}$,
\begin{align}
\notag
\mathrm{V}_0 &= ( A_0 \stackrel{l_0}{\leftarrow} B_0 \stackrel{r_0}{\to} A_1 ) , \\
\mathrm{V}_1 &= ( A_0^\prime \stackrel{l_1}{\leftarrow} B_1 \stackrel{r_1}{\to} A_1^\prime ) . \notag
\end{align}
We define a preorder $\mathrm{V}_1 \triangleright \mathrm{V}_0$ if $A_i = A_i^\prime, ~i = 0,1$ and there exists an epimorphism $e : B_1 \to B_0$ such that $l_0 \circ e = l_1$ and $r_0 \circ e = r_1$.
\end{Defn}

\begin{Defn}
We define a relation $\approx$ of spans in $\mathcal{A}$ by the equivalence relation generated by the preorder $\triangleright$.
\end{Defn}

\begin{remark}
\label{202201261209}
The equivalence relation $\approx$ is streamlined as follows.
Let $\mathrm{V}_0, \mathrm{V}_1$ be spans in $\mathcal{A}$.
Then the following three conditions are equivalent with each other.
\begin{enumerate}
\item
$\mathrm{V}_0 \approx \mathrm{V}_1$.
\item
There exists a lower bound of $\{ \mathrm{V}_0 , \mathrm{V}_1 \}$ with respect to $\triangleright$.
\item
There exists an upper bound of $\{ \mathrm{V}_0 , \mathrm{V}_1 \}$ with respect to $\triangleright$.
\end{enumerate}
Here we sketch a proof of the equivalence of the second and third parts.
Let $\mathrm{V}_0 = ( A_0 \stackrel{l_0}{\leftarrow} B_0 \stackrel{r_0}{\to} A_1 )$, $\mathrm{V}_1 = ( A_0 \stackrel{l_1}{\leftarrow} B_1 \stackrel{r_1}{\to} A_1 ) $ be spans in $\mathcal{A}$.
Suppose that a span $\mathrm{V}_2 = ( A_0 \stackrel{l_2}{\leftarrow} B_2 \stackrel{r_2}{\to} A_1 )$ is a lower bound of $\{ \mathrm{V}_0 , \mathrm{V}_1 \}$.
Let $e : B_0 \to B_2$ and $e^\prime : B_1 \to B_2$ be epimorphisms giving $\mathrm{V}_0 \triangleright \mathrm{V}_2$ and $\mathrm{V}_1 \triangleright \mathrm{V}_2$ respectively.
Let $f : B_0 \oplus B_1 \to B_2$ be the addition of $e$ and $e^\prime$.
We set $g : B_3 \to B_0 \oplus B_1$  by the canonical monomorphism associated with a kernel object $B_3 = \mathrm{Ker} (f)$, and $p_i : B_0 \oplus B_1 \to B_i$ by the $i$-th projection.
Then the span $\mathrm{V}_3 = (A_0 \stackrel{l_0 \circ p_0 \circ g}{\leftarrow} B_3 \stackrel{r_0 \circ p_1\circ g}{\to} A_1)$ is an upper bound of $\{ \mathrm{V}_0 , \mathrm{V}_1 \}$.
Indeed, $p_0 \circ g$ and $(- p_1 \circ g)$ are epimorphisms which give rise to $\mathrm{V}_3 \triangleright \mathrm{V}_0, \mathrm{V}_1$.
Thus, (2) implies (3).
We leave the proof of (2) starting from (3) to readers.
It is obvious that both of (2) and (3) imply (1).
Finally, (1) implies both of (2) and (3) due to the equivalence of (2) and (3).
\end{remark}

\begin{Example}
Let $A$ be an object of $\mathcal{A}$.
It induces a trivial span $t (A) = (0 \leftarrow A \to 0 )$.
Then we have $t (A) \triangleright t (0)$.
Hence, $t (A) \approx t (0) \approx t (B)$ for any object $A,B$.
In other words, any span from a zero object to itself is unique up to $\approx$.
\end{Example}

\begin{Example}
We investigate a span which is slightly nontrivial up to $\approx$.
Let $f : A \to C$ and $g: B \to C$ be any monomorphisms in $\mathcal{A}$.
Consider spans $\mathrm{V}_0 = ( 0 \leftarrow A \stackrel{f}{\to} C )$ and $\mathrm{V}_1 = ( 0 \leftarrow B \stackrel{g}{\to} C )$.
Then $\mathrm{V}_0 \approx \mathrm{V}_1$ if and only if there exists an isomorphism $h : A \to B$ such that $g \circ h = f$.
In fact, $\mathrm{V}_0 \approx \mathrm{V}_1$ implies that there exists an upper bound of $\{ \mathrm{V}_0 , \mathrm{V}_1 \}$ by Remark \ref{202201261209}.
Let $\mathrm{V}_2 = ( 0 \leftarrow X \stackrel{k}{\to} C)$ be an upper bound.
There exists an epimorphism $e : X \to A$ which gives $\mathrm{V}_2 \triangleright \mathrm{V}_0$, i.e. $f \circ e = k$.
Note that $f, e$ give a factorization of $k$ with an image $f$ and a coimage $e$.
Similarly we have a factorization of $k$ via $g$.
Then the claim follows from the uniqueness of factorizations.
\end{Example}

\begin{Defn}
\label{202202031351}
We define a category $\mathsf{Sp}^\approx ( \mathcal{A} )$.
The objects are arbitrary objects in $\mathcal{A}$.
For objects $A,B$, the morphisms from $A$ to $B$ are $\approx$-equivalence classes of spans from $A$ to $B$, i.e. $( A \stackrel{f}{\leftarrow} C \stackrel{g}{\to} B )$ for some morphisms $f,g$ in $\mathcal{A}$.
For each object $A$, the identity of $A$ is the $\approx$-equivalence class of $(A \stackrel{1_A}{\leftarrow} A \stackrel{1_A}{\to} A )$.
For spans $\mathrm{V}_0 = ( A_0 \stackrel{f_0}{\leftarrow} B_0 \stackrel{f_1}{\to} A_1 )$ and $\mathrm{V}_1 = ( A_1 \stackrel{g_0}{\leftarrow} B_1 \stackrel{g_1}{\to} A_2 )$, we define a span $\mathrm{V}_1 \circ \mathrm{V}_0 = ( A_0 \stackrel{f_0 \circ h_0}{\leftarrow} C \stackrel{- g_1 \circ h_1}{\to} A_2 )$ where $C$ is a kernel of $(f_1 + g_1) = \nabla_{A_1} \circ (f_1 \oplus g_1) : B_0 \oplus B_1 \to A_1$ and $h_i : C \to B_i$ are the components of the kernel morphism $\mathrm{ker} (f_1 + g_1) : C \to B_0 \oplus B_1$.
Then the composition in $\mathsf{Sp}^\approx ( \mathcal{A} )$ is defined as $[ \mathrm{V}_1 \circ \mathrm{V}_0 ] \stackrel{\mathrm{def.}}{=} [ \mathrm{V}_1 ] \circ [ \mathrm{V}_0 ]$.
Note that it turns out to be well-defined and satisfies the composition axioms.
The category $\mathsf{Sp}^\approx ( \mathcal{A} )$ has a dagger operation induced by $( A \stackrel{f}{\leftarrow} C \stackrel{g}{\to} B )^\dagger = ( B \stackrel{g}{\leftarrow} C \stackrel{f}{\to} A )$.
An arbitrary morphism $f : A \to B$ in $\mathcal{A}$ induces a span $(A \stackrel{1_A}{\leftarrow} A \stackrel{f}{\to} B )$ in $\mathcal{A}$.
The assignment induces a functor $\iota_\mathcal{A}^+ : \mathcal{A} \to \mathsf{Sp}^\approx ( \mathcal{A} )$.
We set $\iota_\mathcal{A}^- = \dagger \circ \iota_\mathcal{A}^+ : \mathcal{A}^\mathrm{op} \to \mathsf{Sp}^\approx ( \mathcal{A} )$.
\end{Defn}

\begin{remark}
\label{202201301542}
Note that a span $( A \stackrel{f}{\leftarrow} C \stackrel{g}{\to} B )$ induces a relation on $A$ and $B$.
In fact, we choose a epi-mono factorization of the induced morphism $h : C \to A \oplus B$ where $h = (f \oplus g) \circ \Delta_C$.
Let $e : C \to R$ and $m : R \to A \oplus B$ be the facorization.
Then the monomorphism $m : R \to A \oplus B$ gives a relation on $A$ and $B$.
The assignment gives a rise to an isomorphism of categories $R : \mathsf{Sp}^\approx ( \mathcal{A} ) \to \mathsf{Rel} ( \mathcal{A} )$.
\end{remark}

\begin{remark}
\label{202202031403}
Note that the relation category $\mathsf{Rel} ( \mathcal{A} )$, especially $\mathsf{Sp}^\approx ( \mathcal{A} )$ and $\mathsf{Cosp}^\approx ( \mathcal{A} )$ below, are locally small.
Indeed, the morphisms from $A$ to $B$ in $\mathsf{Rel} ( \mathcal{A} )$ coincides with subobjects of $A \oplus B$ which form a set by the assumption.
\end{remark}

\begin{Defn}
There is a parallel discussion about {\it cospans} instead of spans in this subsection.
With a slight abuse of notations, we introduce a preorder $\triangleright$ as follows.
For cospans $\Lambda_0 = ( A_0 \stackrel{l_0}{\to} B_0 \stackrel{r_0}{\leftarrow} A_1 )$ and $\Lambda_1 = ( A_0^\prime \stackrel{l_1}{\to} B_1 \stackrel{r_1}{\leftarrow} A_1^\prime )$, we define $\Lambda_1 \triangleright \Lambda_0$ if $A_i = A_i^\prime, ~ i = 0,1$ and there exists a monomorphism $m : B_0 \to B_1$ such that $m \circ l_0 = l_1$ and $m \circ r_0 = r_1$.
Moreover, a discussion as Definition \ref{202202031351} yields a category $\mathsf{Cosp}^\approx ( \mathcal{A} )$.
Denote by $\iota_\mathcal{A}^+ : \mathcal{A} \to \mathsf{Cosp}^\approx ( \mathcal{A} )$ and $\iota_\mathcal{A}^- : \mathcal{A}^\mathrm{op} \to \mathsf{Cosp}^\approx ( \mathcal{A} )$ the associated canonical functors.
\end{Defn}

\begin{remark}
The categories $\mathsf{Cosp}^\approx ( \mathcal{A} )$ and $\mathsf{Sp}^\approx ( \mathcal{A}^\mathrm{op} )$ coincides with each other by definitions.
\end{remark}

Before we close this subsection, we give an isomorphism between $\mathsf{Cosp}^\approx ( \mathcal{A} )$ and $\mathsf{Sp}^\approx ( \mathcal{A} )$.

\begin{Defn}
Let $\Lambda =   ( A_0 \stackrel{f_0}{\to} B \stackrel{f_1}{\leftarrow} A_1 )$ be a cospan in $\mathcal{A}$.
Let $w_\Lambda = \nabla_B \circ (f_0 \oplus f_1 ) : A_0 \oplus A_1 \to B$ where $\nabla_B : B \oplus B \to B$ is the codiagonal morphism.
Let $g_i : \mathrm{Ker} ( w_\Lambda ) \to A_i$ be the composition of the projection $A_0 \oplus A_1 \to A_i$ and the kernel morphism $\mathrm{Ker} ( w_\Lambda ) \to A_0 \oplus A_1$.
We define the {\it transposition} $\mathrm{T}( \Lambda )$ by a span
\begin{align}
\notag
\mathrm{T}( \Lambda ) \stackrel{\mathrm{def.}}{=} \left( A_0 \stackrel{g_0}{\leftarrow} \mathrm{Ker} ( w_\Lambda ) \stackrel{g_1}{\to} A_1 \right) .
\end{align}
We define a dual notion called the {\it transposition} by a induced cospan $\mathrm{T} (\mathrm{V} )$ for a span $\mathrm{V}$ in $\mathcal{A}$ with a slight abuse of notations.
It turns out that the transpositions induce functors from $\mathsf{Sp}^\approx ( \mathcal{A} )$ to $\mathsf{Cosp}^\approx ( \mathcal{A} )$ and vice versa.
The functors are inverse to each other by the following lemma.
\end{Defn}

\begin{Lemma}
For a cospan $\Lambda$, we have $\mathrm{T} \mathrm{T} ( \Lambda ) \triangleleft \Lambda$.
Especially, we have $\mathrm{T} \mathrm{T} ( \Lambda ) \approx \Lambda$.
\end{Lemma}
\begin{proof}
Let $\Lambda =   ( A_0 \stackrel{f_0}{\to} B \stackrel{f_1}{\leftarrow} A_1 )$.
For the transposition $\mathrm{T}( \Lambda )$ as above, we set $v_{\mathrm{T} ( \Lambda )} = (g_0 \oplus g_1) \circ \Delta$ which is nothing but the kernel morphism of $w_\Lambda$.
Hence, we obtain a morphism $s :  \mathrm{Cok} ( v_{\mathrm{T} ( \Lambda )} )\to B$ with the commutative diagram (\ref{202112031217}).
By the exactness of the upper row, the morphism $s$ is a monomorphism.
Then $s$ gives $\mathrm{T} \mathrm{T} ( \Lambda ) \triangleleft_\mathcal{T} \Lambda$.
\begin{equation}
\label{202112031217}
\begin{tikzcd}[sep=scriptsize]
\mathrm{Ker} ( w_\Lambda ) \ar[r, "v_{\mathrm{T} ( \Lambda )}"] \ar[d, equal] & A_0 \oplus A_1 \ar[r, "w_\Lambda"] \ar[d, equal] & B \\
\mathrm{Ker} ( w_\Lambda ) \ar[r] & A_0 \oplus A_1 \ar[r, "w_{\mathrm{T}\mathrm{T}(\Lambda)}"] & \mathrm{Cok} ( v_{\mathrm{T} ( \Lambda )} ) \ar[u, "s"]
\end{tikzcd}
\end{equation}
\end{proof}

\begin{remark}
Recall Definition \ref{202201271145}.
One can deduce that some common structures are used to construct $\mathsf{Sp}^\approx ( \mathcal{A} )$, $\mathsf{Cosp}^\approx ( \mathcal{A} )$ and the cospan categories of spaces.
Furthermore, both of $\mathsf{Sp}^\approx ( \mathcal{A} )$ and $\mathsf{Cosp}^\approx ( \mathcal{A} )$ are characterized as a universal PN-diagram as an analogue of Theorem \ref{202201172148}.
The phenomena are ubiquitous in algebraic topology and will be dealt with in the sequel.
\end{remark}

\subsection{Cospanical and spanical extensions}
\label{202201301629}

In this subsection, we prove our main theorem by Theorem \ref{202201172145}.
The main theorem is reformulated by using following terminologies.

\begin{Defn}
\label{201912240846}
Let $H : \mathcal{H}_d \to \mathcal{A}$ be a functor.
\begin{enumerate}
\item
A {\it cospanical extension} of $H$ is a dagger-preserving  functor $F : \mathsf{Cosp}^{\simeq}_{\leq (d+1)} (\mathsf{CW}^\mathsf{fin}_{\ast} ) \to \mathsf{Cosp}^{\approx} ( \mathcal{A})$ with the following commutative diagram.
\begin{equation}
\notag
\begin{tikzcd}[sep=scriptsize]
\mathcal{H}_d
\ar[r, "H"]
\ar[d, "\iota_d^+"]
&
\mathcal{A}
\ar[d, "\iota_\mathcal{A}^+"]
\\
\mathsf{Cosp}^{\simeq}_{\leq (d+1)} (\mathsf{CW}^\mathsf{fin}_{\ast} )
\ar[r, "F"]
&
\mathsf{Cosp}^{\approx} ( \mathcal{A})
\end{tikzcd}
\end{equation}
\item
A {\it spanical extension} of $H$ is a dagger-preserving functor $F : \mathsf{Cosp}^{\simeq}_{\leq (d+1)} (\mathsf{CW}^\mathsf{fin}_{\ast} ) \to \mathrm{Sp}^{\approx} ( \mathcal{A})$ with the following commutative diagram.
\begin{equation}
\notag
\begin{tikzcd}[sep=scriptsize]
\mathcal{H}_d
\ar[r, "H"]
\ar[d, "\iota_d^+"]
&
\mathcal{A}
\ar[d, "\iota_\mathcal{A}^+"]
\\
\mathsf{Cosp}^{\simeq}_{\leq (d+1)} (\mathsf{CW}^\mathsf{fin}_{\ast} )
\ar[r, "F"]
&
\mathrm{Sp}^{\approx} ( \mathcal{A})
\end{tikzcd}
\end{equation}
\end{enumerate}
If $H$ is a symmetric monoidal functor, then {\it symmetric monoidal (co)spanical extensions} are defined by dagger-preserving symmetric monoidal functors satsifying the above axioms.
\end{Defn}

\begin{Example}
In DWFQ model based on a finite abelian group $G$, the construction of TQFT's is deeply related with the following cospan induced by a cobordism $N$ from $M_0$ to $M_1$ :
\begin{align}
\notag
\mathbb{C}^{H^1(M_0;G)} \to \mathbb{C}^{H^1(N;G)} \leftarrow \mathbb{C}^{H^1(M_1;G)} .
\end{align}
This assignment yields the cospanical extension of the Mayer-Vietoris functor $\mathbb{C}^{H^1(-;G)}$.
On the one hand, the following span naturally arises from the cobordism where $\partial N \cong M_0 \amalg M_1$ is the boundary of $N$ :
\begin{align}
\notag
\mathbb{C}^{H^1(M_0;G)} \leftarrow \mathbb{C}^{H^2(N , \partial N ;G)} \to \mathbb{C}^{H^1(M_1;G)} 
\end{align}
The assignment of the induced spans yields the spanical extension.
\end{Example}

\begin{Corollary}
\label{202201271352}
For an extended natural number $d = 0, 1, \cdots \infty$, let $E$ be an $\mathcal{A}$-valued MV functor with $\dim \leq (d+1)$.
\begin{enumerate}
\item
There exists a unique cospanical extension $\hat{E}$ of $i_d^\ast (E) = E \circ i_d$.
\item
There exists a unique spanical extension $\check{E}$ of $\Sigma_d^\ast (E) = E \circ \Sigma_d$.
\end{enumerate}
Furthermore, if $d \geq 1$, then there exist a spanical extension of $\Sigma_d^\ast (E)$ and a cospanical extension of $i_{d-1}^\ast ( \Sigma_d^\ast ( E ))$, and we have the following :
\begin{align}
\notag
\mathrm{T} \circ j_d^\ast ( \check{E} ) = \hat{\left(\Sigma_d^\ast ( E )\right)} .
\end{align}
\end{Corollary}
\begin{remark}
\label{202201301549}
If $E$ is pointed (see Remark \ref{202201301359}), then there exist symmetric monoidal (co)spanical extensions subject to the above conditions.
It is proved in a parallel way.
Recall the isomorphism $R : \mathsf{Sp}^\approx (\mathcal{A} ) \to \mathsf{Rel} ( \mathcal{A} )$ in Remark \ref{202201301542}.
The composition $R \circ \mathrm{T}$ gives an isomorphism $\mathsf{Cosp}^\approx (\mathcal{A} ) \to \mathsf{Rel} ( \mathcal{A} )$.
The main theorem is immediate from Corollary \ref{202201271352} up to these isomorphisms.
\end{remark}
\begin{proof}[Proof of Corollary \ref{202201271352}]
Let $F^+ = \iota_\mathcal{A}^+ \circ E : \mathcal{H}_{d+1} \to \mathsf{Cosp}^{\approx} ( \mathcal{A})$ and $F-  = \iota_\mathcal{A}^- \circ E^\mathrm{op} : \mathcal{H}_{d+1}^\mathrm{op} \to \mathsf{Cosp}^{\approx} ( \mathcal{A})$.
They induce a MV PN diagram in $\mathsf{Cosp}^{\approx} ( \mathcal{A})$ with $\dim \leq (d+1)$ :
\begin{align}
\notag
F^+ : \mathcal{H}_{d+1} \longrightarrow \mathsf{Cosp}^{\approx} ( \mathcal{A}) \longleftarrow \mathcal{H}_{d+1}^\mathrm{op} : F^- .
\end{align}
By Theorem \ref{202201172145}, there exists a unification $\hat{E} = i_d^\cup ( F^\pm )$ relative to the embedding functor $i_d$.
Put $\hat{E}^\prime = \dagger \circ \hat{E} \circ \dagger$, i.e. $\hat{E}^\prime ( f ) = \hat{E} (f^\dagger )^\dagger$ for any morphism $f$ in $\mathsf{Cosp}^\simeq_{\leq (d+1)} ( \mathsf{CW}^\mathsf{fin}_\ast )$.
It is easy to verify that $\hat{E}^\prime$ is a unification of $F^\pm$ relative to $i_d$.
By the uniqueness, we have $\hat{E} = \hat{E}^\prime$.
Hence, $\hat{E}$ is a dagger-preserving functor.
The uniqueness of a cospanical extension is immediate from that of the unification.
Similarly the claim with respect to a spanical extension is proved.
The final statement follows from the last claim in Theorem \ref{202201172145}.
\end{proof}

\bibliography{A_cospan_category_of_abelian_category}{}
\bibliographystyle{plain}

\end{document}